\documentclass[reqno,10pt,a4paper]{amsart}
\usepackage[margin=2.7cm]{geometry}
\usepackage[foot]{amsaddr}
\usepackage{cite}
\usepackage{newpxtext}
\usepackage{cabin}% sans serif
\usepackage[varqu,varl]{inconsolata}% sans serif typewriter
\usepackage[bigdelims,vvarbb]{newpxmath}% bb from STIX
\usepackage[cal=cm,bb=esstix,bbscaled=1.05]{mathalfa}% mathcal
\usepackage{booktabs}
\usepackage{verbatimbox}
\usepackage{caption}
\usepackage{lscape}
\usepackage{xcolor}
\usepackage{hyperref}
\hypersetup{
    colorlinks,
    linkcolor={red!50!black},
    citecolor={blue!50!black},
    urlcolor={blue!80!black}
}
\AtBeginDocument{%
\def\MR#1{}
}
\title{Databases of quantum periods for Fano manifolds}
\newcommand{\PP}{\mathbb{P}}
\newcommand{\CC}{\mathbb{C}}
\newcommand{\QQ}{\mathbb{Q}}

\newcommand{\FF}{\mathbb{F}}

\newcommand{\VV}{\mathbb{V}}

\newcommand{\Ghat}{\widehat{G}}
\DeclareMathOperator{\rf}{ramif}
\DeclareMathOperator{\rk}{rank}
\newcommand{\textkey}[1]{`#1'}
\newcommand{\tablekey}[1]{\texttt{#1}}
\begin{document}
%-------------------------------------------------------------------------------
\author[T.\,Coates]{Tom Coates}
\address{Department of Mathematics\\Imperial College London\\180 Queen's Gate\\London\\SW7 2AZ\\UK}
\email{t.coates@imperial.ac.uk}% Coates
\author[A.\,M.\,Kasprzyk]{Alexander Kasprzyk}
\address{School of Mathematical Sciences\\University of Nottingham\\Nottingham\\NG7 2RD\\UK}
\email{a.m.kasprzyk@nottingham.ac.uk}% Kasprzyk
%-------------------------------------------------------------------------------
\keywords{Databases, quantum period, Fano manifolds}
%-------------------------------------------------------------------------------
\begin{abstract}
Fano manifolds are basic building blocks in geometry -- they are, in a precise sense, atomic pieces of shapes. The classification of Fano manifolds is therefore an important problem in geometry, which has been open since the~1930s. One can think of this as building a Periodic Table for shapes. A recent breakthrough in Fano classification involves a technique from theoretical physics called Mirror Symmetry. From this perspective, a Fano manifold is encoded by a sequence of integers: the coefficients of a power series called the regularized quantum period. Progress to date has been hindered by the fact that quantum periods require specialist expertise to compute, and descriptions of known Fano manifolds and their regularized quantum periods are incomplete and scattered in the literature. We describe databases of regularized quantum periods for Fano manifolds in dimensions up to four. The databases in dimensions one, two, and three are complete; the database in dimension four will be updated as new four-dimensional Fano manifolds are discovered and new regularized quantum periods computed. 
\end{abstract}
%-------------------------------------------------------------------------------
\maketitle
%-------------------------------------------------------------------------------
\section*{Background \& summary}
\label{sec:background}
%-------------------------------------------------------------------------------
Algebraic geometry describes geometric shapes as the solution sets of systems of polynomial equations, and lets us exchange geometric understanding of a shape~$X$ with structural understanding of the equations that define it. This geometric analysis has become vital in mathematics and physics, with applications as diverse as motion planning for robots~\cite{HustySchrocker2010,MeniniPossieriTornambe2021}, optimisation~\cite{BlekhermanParriloThomas2013,FotiouRostalskiParriloMorari2006}, algebraic statistics~\cite{PachterSturmfels2005,DrtonSturmfelsSullivant2009}, coding~\cite{vanLintvanderGeer1988,TsfasmanVladut1991}, gauge theory~\cite{AtiyahDrinfeldHitchinManin1978,Hitchin1983}, and string theory~\cite{GreenePlesser1990,CandelasdelaOssaGreenParkes1991,Douglas2001}. The Minimal Model Program decomposes shapes defined by polynomial equations into basic pieces~\cite{Kollar1988,KollarMori1998}. After birational transformations, which are modifications on subsets with zero volume (and codimension at least one), any such shape decomposes into `atomic pieces' of three types: positively curved, flat, and negatively curved. Fano manifolds are the positively curved smooth pieces. These are arguably the most important of the basic pieces, because analysis of basic pieces that are flat~(Calabi\nobreakdash--Yau) or negatively curved~(general type) typically starts by realising them inside a Fano manifold by imposing additional equations. For example, four-dimensional Fano manifolds naturally contain three-dimensional Calabi--Yau manifolds, as their `anticanonical sections', which are of particular importance in constructing models of spacetime in Type~II string theory~\cite{AspinwallGreeneMorrison1994,Polchinski2005}.

There are finitely many Fano manifolds, up to deformation, in each dimension~\cite{KollarMiyaokaMori1992}. There is exactly~1 one-dimensional Fano manifold, the Riemann sphere~$\PP^1$; this follows from the Klein--Poincar\'e Uniformization Theorem. The two-dimensional case is also classical~\cite{DelPezzo1887}; there are~10 deformation families of two-dimensional Fano manifolds. The classification of three-dimensional Fano manifolds, due to Fano~\cite{Fano1947}, Iskovskikh~\cite{Iskovskikh1977,Iskovskikh1978,Iskovskikh1979}, and Mori--Mukai~\cite{MoriMukai1982,MoriMukai1982Erratum} is one of the triumphs of 20th-century algebraic geometry; there are~105 deformation families. The classification of Fano manifolds in dimensions four and higher is still far from understood.

Recently a new approach to the classification of Fano manifolds has been proposed, which relies on a conjectural link between Fano manifolds and Laurent polynomials provided by Mirror Symmetry~\cite{CoatesCortiGalkinGolyshevKasprzyk2013}. This has already led to new classification results~\cite{AkhtarCoatesCortiHeubergerKasprzykOnetoPetracciPrinceTveiten2016} and the discovery of many new four-dimensional Fano manifolds~\cite{CoatesKasprzykPrince2015,Kalashnikov2019, CoatesKasprzykPrince2019}. The key invariant to be analysed is the \emph{regularized quantum period} of a Fano manifold~$X$. This is a power series
\begin{equation}
    \label{eq:Ghat}
    \Ghat_X(t) = 1 + \sum_{d=2}^\infty c_d t^d 
\end{equation}
where~$c_d = r_d d!$ and~$r_d$ is a certain genus-zero Gromov--Witten invariant of~$X$. Intuitively speaking,~$r_d$ is the number of degree-$d$ rational curves in~$X$ that pass through a fixed generic point of~$X$ and have a certain constraint on their complex structure. In general~$r_d$ is a rational number, because curves with a symmetry group of order~$n$ are counted with weight~$1/n$, but in all known cases the coefficients of~$\Ghat_X$ are integers. Gromov--Witten invariants remain constant under deformation of $X$, so~$\Ghat_X$ is a deformation invariant of~$X$. 

This paper describes databases of regularized quantum periods of Fano manifolds in dimensions one, two, three, and~four~\cite{1Ddata,2Ddata,3Ddata,4Ddata}. The databases were prepared by aggregating and standardising existing descriptions of regularized quantum periods in the literature, and computing the regularized quantum periods of various four-dimensional Fano manifolds where geometric constructions were known. The databases in dimensions up to three are complete; the database in dimension four will be updated as further four-dimensional Fano manifolds are discovered and new regularized quantum periods computed. This data will be useful to mathematicians and physicists who wish to identify Fano manifolds, to construct Calabi--Yau manifolds, or to investigate further the connection between Fano manifolds and Mirror Symmetry.
%-------------------------------------------------------------------------------
\section*{Methods}
\label{sec:methods}
%-------------------------------------------------------------------------------
Let~$X$ be a Fano manifold. General theory implies that there is a differential operator
\begin{equation}
    \label{eq:PF}
    L = \sum_{k=1}^N l_k t^{n_k} D^{m_k}
\end{equation}
where the~$l_k$ are non-zero rational numbers and~$D = t \frac{d}{dt}$, such that~$L \Ghat_X\equiv 0$; see e.g.~\cite[Theorem~4.3]{CoatesCortiGalkinGolyshevKasprzyk2013}. We normalise the operator~$L$ in~\eqref{eq:PF} by taking the~$l_k$ to be integers with no common factor and insisting that~$l_M > 0$, where~$M$ is the index~$k$~such that the pair~$(m_k, n_k)$ is lexicographically maximal. The identity~$L \Ghat_X\equiv 0$ translates into a recurrence relation for the coefficients~$c_d$ in~\eqref{eq:Ghat} which determines all of the~$c_d$ from the first few of them, and typically just from~$c_0$ and~$c_1$; here and henceforth we write~$c_0 = 1$ and~$c_1 = 0$. 

Each entry in our databases describes the regularized quantum period of a Fano manifold. Two regularized quantum periods are assumed to be distinct unless it has been proven that they agree (i.e.~that all the infinitely many coefficients~$c_d$ in~\eqref{eq:Ghat} agree). For example, two period sequences are the same if they have been shown to have the same differential operator, or if the Fano manifolds they are derived from have been shown to be deformation equivalent. More precisely, an entry in one of the databases corresponds to a sequence~$c_0, c_1, \ldots$ of coefficients of the regularized quantum period~\eqref{eq:Ghat}, where in some cases only finitely many terms~$c_d$ may be known, along with information about the geometric origins of this period sequence. In dimensions one, two, and three the classification of Fano manifolds is known, and each deformation class of Fano manifolds has a distinct regularized quantum period. Thus in dimensions one, two, and three we can also interpret each entry in our databases as corresponding to a deformation class of Fano manifolds. In dimension four the classification of Fano manifolds is unknown, and it is possible that there exist Fano manifolds~$X_1$ and~$X_2$ that have the same regularized quantum period but that are not deformation equivalent. (No such examples are known, in any dimension, and we expect that no such examples exist.) One should therefore take care to think of each entry in the four-dimensional database as representing \emph{a regularized quantum period sequence of a Fano manifold} rather than a deformation class of Fano manifolds. It may happen that two period sequences in the four-dimensional database are later proven to be equal. When this happens, we will update the database as described in the `Data Records' section.

The following methods were used to compute regularized quantum periods and the differential operators that annihilate~them.

\subsubsection*{Mirror symmetry and the Lairez algorithm}
For~$X$ a smooth Fano toric variety or toric complete intersection, mirror constructions by Givental~\cite{Givental1998} and Hori--Vafa~\cite{HoriVafa2000} give a Laurent polynomial~$f$ that corresponds to~$X$ under Mirror Symmetry -- see e.g.~\cite[\S5]{CoatesKasprzykPrince2015} for a summary of this. Given such a Laurent polynomial, one can compute a differential operator~$L$ such that~$L \Ghat_X\equiv 0$ using Lairez's generalised Griffiths--Dwork algorithm~\cite{Lairez2016}. These steps are implemented in the Fanosearch software library~\cite{magma-core}; see the `Code availability' section.

\subsubsection*{Products}
If the regularized quantum periods of Fano manifolds~$X_1$ and~$X_2$ are known then the regularized quantum period of the product~$X_1 \times X_2$ is determined by~\cite[Corollary~E.4]{CoatesCortiGalkinKasprzyk2016}. Furthermore if~$X_1$ and~$X_2$ correspond under Mirror Symmetry to, respectively, the Laurent polynomials~$f_1 \in \CC[x_1^{\pm 1},\ldots,x_k^{\pm 1}]$ and~$f_2 \in \CC[y_1^{\pm 1},\ldots,y_\ell^{\pm 1}]$, then the product~$X_1 \times X_2$ corresponds under Mirror Symmetry to~$f_1 + f_2 \in \CC[x_1^{\pm 1},\ldots,x_k^{\pm 1},y_1^{\pm 1},\ldots,y_\ell^{\pm 1}]$. One can then use the Lairez algorithm to compute a differential operator~$L$ such that~$L \Ghat_{X_1 \times X_2}\equiv 0$.

\subsubsection*{Numerical linear algebra}

Given sufficiently many coefficients~$c_0,\ldots,c_M$ in~\eqref{eq:Ghat}, one can find a recurrence relation satisfied by this sequence using linear algebra. This determines a candidate for the differential operator~$L$ in~\eqref{eq:PF}. We applied this method to the Strangeway fourfolds, where closed formulas for the~$c_i$ are known~\cite{CoatesGalkinKasprzykStrangeway2014}. When the number~$M$ of coefficients involved is large, the linear system that determines the recurrence relation becomes highly overdetermined and so we can be confident that the operator~$L$ is correct. 

\subsubsection*{A note on rigour}

When computing differential operators~\eqref{eq:PF}, neither the approach based on the Lairez algorithm nor the approach based on numerical linear algebra provides a proof that~$L \Ghat_X\equiv 0$. In the latter case this is because the method cannot do so; in the former case this is because of an implementation detail in Lairez's algorithm: certain calculations over~$\QQ(t)$ are made using reductions to~$\FF_p(t)$ for randomly-chosen primes~$p$ followed by a reconstruction step, and there is a (very small) probability of erroneous reconstruction. In each case, however, the probability of error is tiny. 

In more detail: the Lairez algorithm computes the differential operator~$L$ by first computing a certain connection matrix~$M$ with entries in~$\QQ(t)$. This matrix determines~$L$ uniquely. In the implementation of the algorithm that we used, the matrix~$M$ is reconstructed from its reductions~$M_i$ to~$\FF_{p_i}(t)$ for a sequence~$p_1,\ldots,p_k$ of randomly-chosen 32-bit primes~$p_i$. The entries of~$M$ are recovered from the entries of~$M_i$ using rational reconstruction~\cite{KG83}, where we increase the number~$k$ of primes until the reconstruction stabilises. Let us consider a heuristic estimate of the probability of a single coefficient~$q \in \QQ$ in a single entry of the matrix~$M$ being erroneously reconstructed from its mod-$p_i$ reductions. Set~$m = p_1 p_2 \cdots p_{k-1}$. We have
\begin{align*}
    q \equiv A \mod m && 
    q \equiv B \mod m p_k && \text{for some~$A$,~$B$}
\end{align*}
and since the reconstruction stabilises it follows that~$B \equiv A + r m \mod m p_k$ for some integer~$r$. There are~$p_k$ possible choices for~$B$ given~$A$ and (since~$p_k$ was chosen uniformly at random among 32-bit primes) it seems reasonable that, if the reconstruction at step~$k-1$ was erroneous, then all possibilities for~$B$ are equally likely. Thus the probability of erroneous stabilisation is~$1/p_k$. This is just one coefficient among many. Assuming that erroneous reconstructions of individual coefficients are independent, and using worst-case sizes for~$M$ and degrees of entries in~$M$ gives a probability of erroneous reconstruction on the order of~$10^{-6}$. Furthermore the operators that we found satisfy a number of stringent checks described in the `Technical Validation' section. It is reasonable to conclude that they are correct. 

\subsection*{Dimension one}

The Fano manifold~$\PP^1$ is toric, and corresponds under Mirror Symmetry to the Laurent polynomial~$x + x^{-1}$. The database of regularized quantum periods for one-dimensional Fano manifolds~\cite{1Ddata}, which contains one record, was constructed from this Laurent polynomial using Lairez's algorithm. To cross-check, one can use Givental's mirror theorem~\cite{Givental1996} to compute the regularized quantum period, finding
\[
    \Ghat_{\PP^1} = \sum_{d=0}^\infty t^{2d} \frac{(2d)!}{d!d!}
\]
It is then elementary to check that~$L \Ghat_{\PP^1}\equiv 0$, where~$L = (4t^2-1)D + 4t^2$.

\subsection*{Dimension two}

The database in dimension two~\cite{2Ddata} was constructed by applying the Givental/Hori--Vafa  mirror construction and the Lairez algorithm to the models of two-dimensional Fano manifolds as toric complete intersections given in~\cite[\S G]{CoatesCortiGalkinKasprzyk2016}.

\subsection*{Dimension three}

Regularized quantum periods for three-dimensional Fano manifolds are known, as are Laurent polynomials that correspond to each three-dimensional Fano manifold under Mirror Symmetry~\cite{CoatesCortiGalkinKasprzyk2016}. For~89 of the~105 deformation families, this correspondence follows from the Givental/Hori--Vafa construction. In the remaining cases, which are those in Table~1 of~\cite[Appendix~A]{CoatesCortiGalkinKasprzyk2016} where `Method' is equal to `Abelian/non-Abelian correspondence' or `Quantum Lefschetz with mirror map', the correspondence is conjectural but is supported by strong numerical evidence, including the computation of the first several hundred terms of the expansion~\eqref{eq:Ghat}. The database in dimension three~\cite{3Ddata} was constructed by applying Lairez's algorithm to these Laurent polynomials. 

\subsection*{Dimension four}
The database of regularized quantum periods for four-dimensional Fano manifolds~\cite{4Ddata} was constructed as follows.

\subsubsection*{Four-dimensional Fano toric complete intersections} Four-dimensional Fano manifolds that are complete intersections in smooth Fano toric varieties of dimension up to~8 have been classified~\cite{CoatesKasprzykPrince2015}. The electronic supplementary material for that paper also provides the regularized quantum periods and differential operators~\eqref{eq:PF} in machine-readable form. These differential operators were normalised as discussed after equation~\eqref{eq:PF} and added to the database.

\subsubsection*{Four-dimensional Fano manifolds with classical constructions}
Reference~\citen{CoatesGalkinKasprzykStrangeway2014} computes regularized quantum period sequences and differential operators~\eqref{eq:PF} for many four-dimensional Fano manifolds with classical constructions, including all four-dimensional Fano manifolds with Fano index greater than one. These differential operators were normalised according to our conventions and added to the database. That paper also gives regularized quantum period sequences, but not differential operators, for a number of other four-dimensional Fano manifolds:
\begin{itemize}
    \item four-dimensional Fano toric varieties;
    \item products of lower-dimensional Fano manifolds; 
    \item the Strangeway fourfolds.
\end{itemize}
Toric varieties are toric complete intersections, so four-dimensional Fano toric varieties were already included in the database. Differential operators for products of lower dimensional Fano manifolds were computed using the methods discussed in the section `Products' above. Differential operators for the Strangeway fourfolds were computed as discussed in the section `Numerical linear algebra' above. These differential operators were then normalised as discussed after equation~\eqref{eq:PF} and added to the database.

\subsubsection*{Four-dimensional quiver flag zero loci}
Four-dimensional Fano manifolds that are quiver flag zero loci in Fano quiver flag varieties of dimension up to~8 were classified by Kalashnikov~\cite{Kalashnikov2019}, who also computed the coefficients~$c_d$ in~\eqref{eq:Ghat} for these for~$d \leq 15$.  Kalashnikov partitions the four-dimensional Fano quiver flag zero loci into equivalence classes depending on the values of the coefficients~$(c_0,\ldots,c_{15})$~and finds~749~equivalence classes. It is not known, however, whether quiver flag zero loci with the same~$(c_0,\ldots,c_{15})$ are deformation equivalent, so we regard the regularized quantum period sequence represented by~$(c_0,\ldots,c_{15})$ as coming from the four-dimensional Fano manifold specified as the representative of the equivalence class in~\cite[Appendix~B, Table~1]{Kalashnikov2019}. These representatives were added to the database.
%-------------------------------------------------------------------------------
\section*{Data records}
\label{sec:data_records}
%-------------------------------------------------------------------------------
We provide four key-value databases, containing regularized quantum periods for Fano manifolds in dimensions one, two, three, and~four. These databases have been committed to the public domain using a~CC0 license~\cite{CC0}. They are available from the open access data repository Zenodo~\cite{1Ddata,2Ddata,3Ddata,4Ddata}. Zenodo provides versioned~DOIs, and new versions of the databases will be produced if the data needs to be updated -- for example because new four-dimensional quantum periods are computed. Each database is presented as an ASCII text file, called `smooth\_fano\_$N$.txt' where $N$ is the appropriate dimension. Each line of a record in that file contains a key, followed by `: ', followed by a value. Records are separated by a single blank line.

Table~\ref{tab:keys} describes the keys provided by the records in each database and their corresponding values, where we think of each database as a single key-value table. Keys are strings of lower-case characters. Each entry in each database has an~ID, which is a positive integer. IDs are sequential and start from~1, but carry no meaning or information other than to identify that particular entry in that particular database. Each entry in each database also specifies a non-empty sequence of names. These names identify Fano manifolds with that regularized quantum period in various published~(partial) classifications, and hence determine constructions of these Fano manifolds. Names of the form `\texttt{S1 x S2}', where \texttt{S1} and \texttt{S2} are names of Fano manifolds~$X_1$~and~$X_2$, indicate that the corresponding Fano manifold is a product~$X_1 \times X_2$. Fano manifolds can have many different names. The databases of regularized quantum period sequences for one-, two-, and three-dimensional Fano manifolds use names as in Table~\ref{tab:low_dimensional_names}; the database of regularized quantum period sequences for four-dimensional Fano manifolds uses names as in Table~\ref{tab:four_dimensional_names}, and as specified in the file~`README.txt' in the Zenodo dataset~\cite{4Ddata}. 

The differential operator~$L$ in~\eqref{eq:PF} is recorded as parallel sequences of coefficients~$l_k$ and exponents~$(m_k, n_k)$; the corresponding keys in the database are named \textkey{pf\_coefficients} and \textkey{pf\_exponents} to reflect the fact that under Mirror Symmetry~$L$ corresponds to a Picard--Fuchs operator. The databases in dimensions one, two, and~three contain values for all keys except \textkey{duplicate}; in particular this determines the differential operator~$L$ in~\eqref{eq:PF}, and thus all coefficients~$c_d$ in the expansion~\eqref{eq:Ghat} of~$\Ghat_X$. For a number of the entries in the four-dimensional database, the differential operator~\eqref{eq:PF} is unknown; in this case the key-value pairs related to the differential operator are omitted from the database entry.

The databases contain a key \textkey{pf\_proven}, with boolean values, that has the value \texttt{false} for each entry in each database such that~\eqref{eq:PF} is known. This reflects the fact that, as discussed above, at the time of writing it is not proven that the differential operators~\eqref{eq:PF} recorded in the database actually annihilate the corresponding regularized quantum periods. If this situation changes in the future -- for example, if the certificated version of the Lairez algorithm described in~\cite[\S7.3]{Lairez2016} is implemented -- then we will make new versions of the databases available.

The presence of the key \textkey{duplicate} in an entry indicates that this entry~$E$ is a duplicate of another entry in the same database, with the indicated~ID, and that the entry~$E$ will not change further as the database is updated. The key \textkey{duplicate} is not (and will never be) present in any of the databases in dimensions one, two, or three; at the time of writing it is also not present in the database in dimension four, but it will be added in future updates to indicate when two entries in that database have been proven to coincide.
%-------------------------------------------------------------------------------
\section*{Technical validation}
\label{sec:validation}
%-------------------------------------------------------------------------------
The records in our databases satisfy a number of consistency checks. These can be verified, for example, using the computational algebra system Magma~\cite{BosmaCannonPlayoust1997} and the Fanosearch software library~\cite{magma-core}; see Table~\ref{tab:magma_core} for the relevant function names, and the software documentation for these functions for the arguments and parameters required. Firstly, the period sequences in each database are annihilated by the corresponding differential operators~\eqref{eq:PF} whenever they are known. Secondly. the differential operators~$L$~from~\eqref{eq:PF} are expected to be of Fuchsian type, that is, to have only regular singular points. This is an extremely delicate condition on the coefficients~$l_k$, and can be checked by exact computation; in particular calls to \texttt{RamificationData(L)} or~\texttt{RamificationDefect(L)} in the Fanosearch Magma library will raise an error if \texttt{L} is not Fuchsian. All the differential operators in the databases in dimensions one, two, and~three are Fuchsian. Checking this for some of the entries in the four-dimensional database is prohibitively expensive, because it involves linear algebra over the splitting field of the symbol of the differential operator and this symbol can be of very high degree, but all of the four-dimensional operators for which the calculation was possible are Fuchsian, and all of them have regular singularities at those singular points defined over number fields of low degree. 

As a further check, the differential operator~\eqref{eq:PF} is expected to be of \emph{low ramification} in the following sense~\cite{CoatesCortiGalkinGolyshevKasprzyk2013}. Let~$S \subset \PP^1$ be a finite set and~$\VV \to \PP^1 \setminus S$ a local system.  Fix a basepoint~$x \in \PP^1 \setminus S$.  For~$s \in S$, choose a small loop that winds once anticlockwise around~$s$ and connect it to~$x$ via a path, thereby making a loop~$\gamma_s$ about~$s$ based at~$x$.  Let~$T_s \colon \VV_x \to \VV_x$ denote the monodromy of~$\VV$  along~$\gamma_s$.  The \emph{ramification} of~$\VV$ is:
\begin{align*}
\rf(\VV) = \sum_{s \in S} \dim\Big(\VV_x/{\VV_x}^{\!\!\!T_s}\Big) && \text{where~${\VV_x}^{\!\!\!T_s}$ is the~$T_s$-invariant part of~$\VV_x$.}
\end{align*}
The \emph{ramification defect} of~$\VV$ is the quantity~$\rf(\VV) - 2\rk(\VV)$.  Non-trivial irreducible local systems~$\VV \to \PP^1 \setminus S$ have $\rf(\VV) \geq 2\rk(\VV)$, and hence have non-negative ramification defect.  A local system of ramification defect zero is called~\emph{extremal}. The \emph{ramification} (and respectively \emph{ramification defect}) of a differential operator~$L$ is the ramification (and respectively ramification defect) of the local system of solutions~$L f \equiv 0$. All the differential operators in the databases in dimensions one, two, and~three are of low ramification; indeed with the exception of the two-dimensional Fano manifolds with names~`dP(7)' and~`dP(8)', which have ramification defect~1, all of these operators are extremal. Computing the ramification for some of the operators in the four-dimensional database is prohibitively expensive, for the same reason as before, but all of the four-dimensional operators for which the calculation was possible are of low ramification, and many of them extremal or of ramification defect~1.
%-------------------------------------------------------------------------------
\section*{Usage notes}
\label{sec:usage}
%-------------------------------------------------------------------------------
The key-value data files provided in the Zenodo datasets are easy to parse in any computational algebra system, and in particular can be parsed using the function \texttt{KeyValueFileProcess} provided by the Fanosearch Magma library~\cite{magma-core}. The four-dimensional database is available for interactive searching via the Graded Ring Database~\cite{GRDB} at 
\begin{center}
    \url{http://grdb.co.uk/forms/smoothfano4}
\end{center}
and programmatically via the Graded Ring Database~API. For example, a request to the~URL
\begin{center}
    \texttt{http://grdb.co.uk/xml/search.xml?agent=curl\&dataid=smoothfano4\&\phantom{XXXXXX}\\
    \hspace{0.6\textwidth}c4=72\&c5=360\&printlevel=1}
\end{center}
will return~XML-formatted data as follows:
\begin{verbnobox}[\hspace{0.1\textwidth}]
<?xml version="1.0"?>
<!-- Graded Ring Database -->
<results numrows="1">
  <result row="1" printlevel="1">
    <id>32</id>
    <names>CKP(31), Obro(4,31)</names>
    <c2>0</c2>
    <c3>18</c3>
    <c4>72</c4>
    <c5>360</c5>
    <c6>2430</c6>
  </result>
</results>
\end{verbnobox} 
Increasing the value of \texttt{printlevel} in the request, to a maximum of~3, will return more information. Users are encouraged to change the value of \texttt{agent} to something more appropriate for their application. For example, a request to the~URL
\begin{center}
    \texttt{http://grdb.co.uk/xml/search.xml?agent=my\_app\&dataid=smoothfano4\&\phantom{XXXXXX}\\
    \hspace{0.6\textwidth}id=340\&printlevel=2}
\end{center}
will return

\begin{verbnobox}[\hspace{0.1\textwidth}]
<?xml version="1.0"?>
<!-- Graded Ring Database -->
<results numrows="1">
  <result row="1" printlevel="2">
    <id>340</id>
    <names>CKK(262), CKP(332)</names>
    <c2>6</c2>
    <c3>6</c3>
    <c4>114</c4>
    <c5>360</c5>
    <c6>3390</c6>
    <period>[1,0,6,6,...]</period>
    <notes>This period sequence is realised by ...</notes>
  </result>
</results>
\end{verbnobox}
Here `\ldots' indicates that some output has been omitted, for readability. The~XML elements \texttt{c2},~\ldots,~\texttt{c6} of the \texttt{result} element record the coefficients~$c_2,\ldots,c_6$ in \eqref{eq:Ghat}; there are also~XML elements of the \texttt{result} element with the same names as the corresponding keys in Table~\ref{tab:keys}. Elements \texttt{pf\_coefficients}, \texttt{pf\_exponents}, and \texttt{pf\_proven} are included only when \texttt{printlevel=3}.

\subsection*{Calabi--Yau differential operators}
When~$X$ is a four-dimensional Fano manifold of Picard rank one, and in certain other sporadic cases, the differential operator~$L$~in~\eqref{eq:PF} is a Calabi--Yau differential operator~\cite{vanStraten2018} of order~4. Thus some of the regularized quantum periods that we describe also appear in the~AESZ table~\cite{AlmkvistvanEnckevortvanStratenZudilin2005,AESZ} of Calabi--Yau differential operators.
%-------------------------------------------------------------------------------
\section*{Code availability}
\label{sec:code_availability}
%-------------------------------------------------------------------------------
The databases were prepared using the Fanosearch software library~\cite{magma-core}, which is freely available under a~CC0 license. The commit hash in that reference records the precise version of the software used. Table~\ref{tab:magma_core} describes intrinsics (i.e.~functions in the computational algebra system Magma~\cite{BosmaCannonPlayoust1997}) provided by that library that can be used to rebuild the database, or to perform the consistency checks described in the Technical Validation section. Lairez's original implementation of his generalised Griffiths--Dwork algorithm is available from GitHub~\cite{Lairez-periods} under a~CeCILL license.
%-------------------------------------------------------------------------------
\section*{Acknowledgements} 
%-------------------------------------------------------------------------------
TC is funded by ERC Consolidator Grant~682603 and EPSRC Programme Grant~EP/N03189X/1. AK is supported by EPSRC~Fellowship~EP/N022513/1. We thank Alessio Corti and Pieter Belmans for providing useful perspectives and~comments. 
%-------------------------------------------------------------------------------

%-------------------------------------------------------------------------------
\appendix
%-------------------------------------------------------------------------------
\begin{landscape}
\section*{Figures \& tables}
\label{sec:tables}
%-------------------------------------------------------------------------------
\begin{table}[htbp]
    \small
    \centering\begin{tabular}{lll}
        \toprule
        \multicolumn{1}{c}{Key} & \multicolumn{1}{c}{Type and Value} \\
        \midrule
        \tablekey{id} &  a positive integer specifying the~ID of the entry \\ 
        \tablekey{period} & a sequence of non-negative integers giving the first few terms~$c_d$ in \eqref{eq:Ghat} \\
        \tablekey{names} & a sequence of strings specifying names for the Fano manifold: see Tables~\ref{tab:low_dimensional_names} and~\ref{tab:four_dimensional_names} \\
        \tablekey{pf\_coefficients} & a sequence of integers~$l_k$ as in~\eqref{eq:PF} \\
        \tablekey{pf\_exponents} & a sequence of length-two sequences~$[m_k, n_k]$ of positive integers as in \eqref{eq:PF} \\
        \tablekey{pf\_proven} & a boolean, true if and only if the operator~\eqref{eq:PF} is proven to annihilate~$\Ghat_X$ \\
        \tablekey{notes} & a string containing further notes on the Fano manifold \\
        \tablekey{duplicate} & a positive integer specifying the~ID of another entry in the same database \\
        \bottomrule \\
    \end{tabular}
    \vspace{-1em}
    \caption{The keys and values in each of the databases}
    \label{tab:keys}
\end{table}

\begin{table}[htbp]
    \small
    \centering\begin{tabular}{ll}
        \toprule
        \multicolumn{1}{c}{Name} & \multicolumn{1}{c}{Meaning} \\
        \midrule
        \tablekey{P1} & one-dimensional projective space \\
        \tablekey{P2} & two-dimensional projective space \\
        \tablekey{P3} & three-dimensional projective space \\
        \tablekey{Q3} & a quadric hypersurface in~$\PP^4$ \\
        \tablekey{dP(k)} & the del~Pezzo surface of degree~$k$ given by the blow-up of~$\PP^2$ in~$9-k$ points\\
        \tablekey{V(3,k)} & the three-dimensional Fano manifold of Picard rank~$1$, Fano index~$1$, and degree~$k$ \\
        \tablekey{B(3,k)} & the three-dimensional Fano manifold of Picard rank~$1$, Fano index~$2$, and degree~$8k$ \\
        \tablekey{MM(r,k)} & the~$k$th entry in the Mori--Mukai list of three-dimensional Fano manifolds of Picard rank~$r$, ordered as in reference~\citen{CoatesCortiGalkinKasprzyk2016} \\
        \bottomrule \\
    \end{tabular}
    \vspace{-1em}
    \caption{Names of entries in the databases for dimensions one, two, and three.}
    \label{tab:low_dimensional_names}
\end{table}

\begin{table}[htbp]
    \small
    \centering\begin{tabular}{ll}
        \toprule
        \multicolumn{1}{c}{Name} & \multicolumn{1}{c}{Meaning} \\
        \midrule
        \tablekey{P4} & four-dimensional projective space \\
        \tablekey{Q4} & a quadric hypersurface in~$\PP^5$ \\
        \tablekey{FI(4,k)} & the four-dimensional Fano manifold of Fano index~$3$ and degree~$81k$ \\
        \tablekey{V(4,k)} & the four-dimensional Fano manifold of Picard rank~$1$, Fano index~$2$, and degree~$16k$ \\
        \tablekey{MW(4,k)} & the~$k$th entry in the table~\cite[Table~12.7]{IskovskikhProkhorov1999} of four-dimensional Fano manifolds of Fano index~$2$ and Picard rank~$\ge2$ \\
        \tablekey{Obro(4,k)} & the~$k$th four-dimensional Fano toric manifold in {\O}bro's classification~\cite{Obro2007}\\
        \tablekey{Str(k)} & the~$k$th Strangeway manifold in reference~\citen{CoatesGalkinKasprzykStrangeway2014} \\
        \tablekey{CKP(k)} & the~$k$th four-dimensional Fano toric complete intersection in reference~\citen{CoatesKasprzykPrince2015} \\
        \tablekey{CKK(k)} & the~$k$th four-dimensional Fano quiver flag zero locus in Appendix~B of reference~\citen{Kalashnikov2019} \\
        \bottomrule \\
    \end{tabular}
    \vspace{-1em}
    \caption{Names of entries in the database for dimension four.}
    \label{tab:four_dimensional_names}
\end{table}

\begin{table}[htbp]
    \small
    \centering\begin{tabular}{ll}
        \toprule
        \multicolumn{1}{c}{Operation} & \multicolumn{1}{c}{Function name} \\
        \midrule
        Read the key-value data files in the Zenodo datasets & \texttt{KeyValueFileProcess} \\
        Compute regularized quantum periods for toric complete intersections & 
        \texttt{PeriodSequenceForCompleteIntersection} \\
        Compute regularized quantum periods for quiver flag zero loci & 
        \texttt{PeriodSequence} \\
        Compute Picard--Fuchs operators using Lairez's algorithm & 
        \texttt{PicardFuchsOperator} applied to a Laurent polynomial \\
        Compute Picard--Fuchs operators using numerical linear algebra &
        \texttt{PicardFuchsOperator} applied to an integer sequence \\
        Apply a Picard--Fuchs operator to a period sequence &
        \texttt{ApplyPicardFuchsOperator} \\
        Compute ramification data for a Fuchsian differential operator & \texttt{RamificationData} \\
        Compute the ramification defect of a Fuchsian differential operator & \texttt{RamificationDefect} \\
        \bottomrule \\
    \end{tabular}
    \vspace{-1em}
    \caption{Useful Magma functions (also known as intrinsics) in the Fanosearch library}
    \label{tab:magma_core}
\end{table}

\end{landscape}
%-------------------------------------------------------------------------------
\end{document}